\documentclass[12pt, twoside]{article}

\usepackage{amsfonts}

\usepackage{amssymb, amsthm, amsmath}

\usepackage{amstext}

\date{} 

\begin{document} 


\centerline{} 

\centerline{} 

\centerline {\Large{\bf A Remark on Compact  Minimal Surfaces in $S^5$  }} 

\centerline{} 

\centerline {\Large{\bf with Non-Negative Gaussian Curvature}} 

\centerline{} 

\centerline{\bf {Rodrigo Ristow Montes}} 

\centerline{} 

\centerline{Departamento de Matem\'atica}

\centerline{Universidade Federal da Para\'iba}

\centerline{BR-- 58.051-900 ~ Jo\~ao Pessoa, P.B., Brazil }

\centerline{ristow@mat.ufpb.br} 

\newtheorem{Theorem}{\quad Theorem}[section] 

\newtheorem{Remark}[Theorem]{\quad Remark} 

\newtheorem{Definition}[Theorem]{\quad Definition} 

\newtheorem{Corollary}[Theorem]{\quad Corollary} 

\newtheorem{Lemma}[Theorem]{\quad Lemma} 

\begin{abstract} 

  In this paper we classify compact minimal surfaces in $S^5$ with non-negative Gaussian curvature using the notion of a contact angle. 

\end{abstract}

{\bf Mathematics Subject Classification:} 53D10 - 53C42 - 53C21 \\ 
 
{\bf Keywords:} contact angle, holomorphic angle, minimal surfaces, contact distribution.


\section{Introduction}
In  \cite{MV} we introduced the notion of contact angle, that can be
considered as  a new geometric invariant useful to investigate
the geometry of immersed surfaces in $S^{3}$. Geometrically, the contact angle
$(\beta)$ is the complementary angle between the contact
distribution and the tangent space of the surface. Also in
\cite{MV}, we deduced formulas for the Gaussian curvature and the Laplacian
of an immersed minimal surface in $S^3$, and we gave a
characterization of the  Clifford Torus as the only minimal
surface in $S^3$
with constant contact angle.\\
We define $\alpha$ to be the angle given by $ \cos \alpha = \langle
ie_1 , v \rangle \nonumber$, where $e_1$ and $v$ are defined in
section~\ref{sec:section2}. The holomorphic angle $\alpha$ is the analogue of
the K\"ahler
angle introduced by Chern and Wolfson in \cite{CW}.\\
Recently,  in \cite{RMV},  we construct a
family of minimal tori in $S^{5}$ with constant contact and
holomorphic angle. These tori are parametrized by the following circle
equation
\begin{equation}\label{eq:equacaoab}
a^2 + \left(b- \frac{\cos\beta}{1+\sin^2\beta}\right)^2 =  2\frac{\sin^4\beta}{(1+\sin^2\beta)^2},
\end{equation}
where  $a$ and $b$ are given in Section~\ref{parametro} (equation
(\ref{eq: segundap})).
In particular, when $a=0$ in \eqref{eq:equacaoab}, we recover the
examples found by Kenmotsu, in \cite{KK2}.  These examples are
defined for $0 < \beta < \frac{\pi}{2}$. Also, when $b=0$ in
\eqref{eq:equacaoab}, we find a new family of minimal tori in
$S^5$, and these tori are defined for $\frac{\pi}{4} < \beta <
\frac{\pi}{2}$. Also, in \cite{RMV}, when $\beta = \frac{\pi}{2}$, we
give an alternative proof of this classification of a Theorem from 
Blair in \cite{Blair}, and Yamaguchi, Kon and Miyahara in
\cite{YKM} for Legendrian minimal surfaces in $S^5$ with
constant Gaussian curvature.\\
In this paper, we will establish a congruence theorem for  minimal surfaces in $S^5$ with constant holomorphic angle ($\alpha$). Using Codazzi-Ricci equations, we prove the following theorem:
\begin{Theorem}\label{teorema1}
Consider $S$ a riemannian surface, $\beta : S \rightarrow ]0,\frac{\pi}{2}[$ a function over $S$ that verifies the following equation:
\begin{eqnarray}
\Delta(\beta) &  =   & \cot\beta ((1+\csc^2\beta)a^2 -\sin^2\alpha) - \tan\beta(|\nabla\beta|^2+12\cos^2\alpha-\sin^2\alpha) \nonumber
\end{eqnarray}
then there exist one and  only one minimal immersion of $S$ into $S^5$ with constant holomorphic angle such that $\beta$ is the contact angle of this immersion.
\end{Theorem}
As a consequence of the Gauss equation, we will provide a theorem of the classification of compact minimal surfaces in $S^5$ supposing that the  Gaussian Curvature is non-negative ($K \geq 0$), then we have the following result:
\begin{Corollary}\label{corolario1}
 Suppose that the holomorphic  angle ($\alpha$) is constant and suppose that $S$ is a compact minimal surface in $S^5$ with constant principal curvature in the direction $e_3$,ie. $a$ is constant, with non-negative Gaussian curvature ($K \geq 0$), then the contact angle ($\beta$) must be constant. 
\end{Corollary}
\section{Contact Angle for Immersed Surfaces in  $S^5$}\label{sec:section2}
Consider in $\mathbb{C}^3$ the following objects:
\begin{itemize}
\item the Hermitian product: $(z,w)=\sum_{j=0}^2 z^j\bar{w}^j$;
\item the inner product: $\langle z,w \rangle = Re (z,w)$;
\item the unit sphere: $S^5=\big\{z\in\mathbb{C}^3 | (z,z)=1\big\}$;
\item the \emph{Reeb} vector field in $S^5$, given by: $\xi(z)=iz$;
\item the contact distribution in $S^5$, which is orthogonal to $\xi$:
\[\Delta_z=\big\{v\in T_zS^5 | \langle \xi , v \rangle = 0\big\}.\]
\end{itemize}
We observe that $\Delta$ is invariant by the complex structure of $\mathbb{C}^3$.

Let now $S$ be an immersed orientable surface in $S^5$.
\begin{Definition}
The \emph{contact angle} $\beta$ is the complementary angle between the
contact distribution $\Delta$ and the tangent space $TS$ of the
surface.
\end{Definition}
Let $(e_1,e_2)$ be  a local frame of $TS$, where $e_1\in
TS\cap\Delta$. Then $\cos \beta = \langle \xi , e_2 \rangle
$. Finally, let $v$ be the unit vector in the direction of  the orthogonal projection  of $e_2$ on $\Delta$,
defined by the following relation
\begin{eqnarray}\label{eq:campoe2}
e_2 = \sin\beta v + \cos\beta \xi.
\end{eqnarray}
\section{Equations for Gaussian curvature and Laplacian of a minimal surface in
  $S^5$ with Constant Holomorphic Angle $\alpha$}\label{parametro}
In this section, we deduce the equations for the Gaussian curvature  and for the Laplacian
of a minimal surface in $S^5$ in terms of the
contact angle and the holomorphic angle.
Consider the normal vector fields
\begin{eqnarray}\label{eq:camposnormais}
e_3                   & = & i \csc \alpha e_1 - \cot \alpha v \nonumber\\
e_4                   & = & \cot \alpha e_1 + i \csc \alpha v \\
e_5                   & = & \csc \beta \xi - \cot \beta e_2
\nonumber
\end{eqnarray}
where $\beta \neq {{0,\pi}}$ and  $ \alpha \neq{{0,\pi}}$.
We will call $(e_j)_{1 \leq \ j \leq 5}$ an  \emph{adapted frame}.

Using (\ref{eq:campoe2}) and (\ref{eq:camposnormais}), we get
\begin{eqnarray}\label{eq:camposdistribution}\
v  =  \sin \beta e_2  - \cos \beta e_5 , \quad  iv   =  \sin \alpha e_4 - \cos \alpha e_1 \\
\xi  =  \cos \beta e_2 + \sin \beta e_5 \nonumber
\end{eqnarray}
It follows from (\ref{eq:camposnormais}) and
(\ref{eq:camposdistribution}) that
\begin{eqnarray}\label{eq:inormais}
ie_1 & = & \cos \alpha \sin \beta e_2 + \sin\alpha e_3 -\cos\alpha \cos \beta
e_5  \\
ie_2 & = & - \cos \beta z - \cos\alpha \sin\beta e_1  + \sin\alpha \sin
\beta e_4 \nonumber
\end{eqnarray}
Consider now the dual basis $(\theta^j)$ of $(e_j)$.
The connection forms $(\theta_k^j)$ are given by
\begin{eqnarray}
De_j    =   \theta_j^k e_k, \nonumber
\end{eqnarray}
and the second fundamental form  with respect to this frame are given by
\begin{equation}
\begin{array}{lclll}
II^j      & = & \theta_1^j \theta^1  +  \theta_2^j \theta^2; \quad
j=3, ..., 5 \nonumber.
\end{array}
\end{equation}
Using (\ref{eq:inormais}) and differentiating $v$ and $\xi$ on the
surface $S$, we get
\begin{eqnarray}\label{eq:dif}
D\xi & = & -\cos\alpha \sin\beta \theta^2 e_1 + \cos \alpha \sin \beta
 \theta^1 e_2 + \sin \alpha \theta^1 e_3 + \sin\alpha \sin \beta \theta^2 e_4\nonumber\\
          && - \cos \alpha \cos \beta \theta^1 e_5,  \\
Dv   & = & (\sin \beta \theta_2^1 - \cos \beta \theta_5^1)e_1 +
 \cos\beta(d\beta-\theta_5^2)e_2 + ( \sin\beta \theta_2^3 - \cos \beta
 \theta_5^3)e_3 \nonumber \\
          && +  ( \sin\beta \theta_4^2 - \cos \beta \theta_5^4)e_4 + \sin
 \beta(d\beta + \theta_2^5)e_5 \nonumber.
\end{eqnarray}
Differentiating $e_3$, $e_4$ and  $e_5$, we have
\begin{eqnarray}\label{eq:intrin}
\theta_3^1 & = &  -\theta_1^3 \nonumber \\
\theta_3^2 & = &  \phantom{-}\sin \beta \theta_4^1 - \cos \beta \sin \alpha
\theta^1 \nonumber \\
\theta_3^4 & = & \phantom{-} \csc \beta \theta_1^2 - \cot \alpha ( \theta_1^3 + \csc \beta
\theta_2^4) \nonumber \\
\theta_3^5 & = & \phantom{-} \cot \beta \theta_2^3 - \csc \beta \sin \alpha \theta^1 \nonumber \\
\theta_4^1 & = &  - \csc \beta \theta_2^3 + \sin \alpha \cot \beta \theta^1 \nonumber \\
\theta_4^2 & = & - \theta_2^4 \nonumber \\
\theta_4^3 & = & \phantom{-} \csc \beta \theta_2^1 + \cot \alpha ( \theta_1^3 + \csc \beta
\theta_2^4) \nonumber \\
\theta_4^5 & = & \phantom{-} \cot \beta \theta_2^4 - \sin \alpha \theta^2 \nonumber \\
\theta_5^1 & = & -\cos \alpha \theta^2 -  \cot \beta \theta_2^1 \nonumber \\
\theta_5^2 & = & \phantom{-} d\beta + \cos \alpha \theta^1 \\  \label{eq: intrin}
\theta_5^3 & = & -\cot \beta \theta_2^3 + \csc \beta \sin \alpha \theta^1 \nonumber\\
\theta_5^4 & = & -\cot \beta \theta_2^4 + \sin \alpha \theta^2 \nonumber
\end{eqnarray}
The conditions of minimality and of symmetry  are equivalent to the following
equations:
\begin{eqnarray}\label{eq:minimalp}
\theta_1^\lambda \wedge \theta^1 +  \theta_2^\lambda \wedge
\theta^2  = 0 = \theta_1^\lambda \wedge \theta^2 -
\theta_2^\lambda \wedge \theta^1.
\end{eqnarray}
On the surface $S$, we consider
\begin{eqnarray}
\theta_1^3 & = & a\theta^1 + b\theta^2 \nonumber
\end{eqnarray}
It follows from (\ref{eq:minimalp}) that
\begin{eqnarray}\label{eq: segundap}
\theta_1^3 & = & \phantom{-} a\theta^1 + b\theta^2\nonumber\\
\theta_2^3 & = & \phantom{-} b\theta^1 - a\theta^2\nonumber\\
\theta_1^4 & = &  (b \csc\beta - \sin\alpha
\cot\beta)\theta^1 - a \csc \beta \theta^2\nonumber \\
\theta_2^4 & = & \phantom{-}  - a \csc\beta \theta^1 - (b
\csc\beta - \sin\alpha \cot\beta)\theta^2 \\
\theta_1^5   &  =  & \phantom{-} d\beta \circ J  - \cos \alpha \theta^2\nonumber\\
\theta_2^5   &  =  &  - d\beta - \cos \alpha \theta^1\nonumber
\end{eqnarray}
where $J$ is the complex structure  of $S$ is given by $Je_1=e_2$ and $Je_2=-e_1$.
Moreover, the normal connection forms are given by:
\begin{eqnarray}\label{eq:normalconexaobetap}
\theta_3^4 & = & - \sec\beta d\beta \circ J + a \cot\alpha \cot^2\beta \theta^1 \nonumber \\
&&+ ( b  \cot\alpha \cot^2\beta -
\cos\alpha \cot\beta \csc\beta + 2 \sec\beta \cos \alpha) \theta^2 \nonumber \\
\theta_3^5 & = & \phantom{-} (b \cot\beta - \csc\beta \sin\alpha) \theta^1 - a
\cot\beta \theta^2 \\
\theta_4^5 & = & - a \cot\beta \csc\beta \theta^1 + ( -b\csc\beta \cot\beta + \sin\alpha (\cot^2\beta -1))\theta^2,\nonumber
\end{eqnarray}
 We consider diagonalizable the  fundamental second  form in the direction $e_3$. We are interested in order to 
study the case where $b=0$. Therefore, we have:
\begin{eqnarray}
\theta_1^3 & = & a\theta^1 \nonumber
\end{eqnarray}
It follows from (\ref{eq: segundap}) that
\begin{eqnarray}\label{eq: segunda}
\theta_1^3 & = & \phantom{-} a\theta^1 \nonumber\\
\theta_2^3 & = & \phantom{-} - a\theta^2\nonumber\\
\theta_1^4 & = & \phantom{-}  - \sin\alpha
\cot\beta\theta^1 - a \csc \beta \theta^2\nonumber \\
\theta_2^4 & = & \phantom{-} - a \csc\beta \theta^1 +  \sin\alpha \cot\beta \theta^2 \\
\theta_1^5   &  =  & \phantom{-} d\beta \circ J  - \cos \alpha \theta^2\nonumber\\
\theta_2^5   &  =  &  - d\beta - \cos \alpha \theta^1\nonumber
\end{eqnarray}
where $J$ is the complex structure  of $S$ is given by $Je_1=e_2$ and $Je_2=-e_1$.
Moreover, the normal connection forms are given by:
\begin{eqnarray}\label{eq:normalconexaobeta}
\theta_3^4 & = & - \sec\beta d\beta \circ J + a \cot\alpha \cot^2\beta \theta^1 \nonumber \\
&&+ (  - \cos\alpha \cot\beta \csc\beta + 2 \sec\beta \cos \alpha) \theta^2 \nonumber \\
\theta_3^5 & = & \phantom{-} - \csc\beta \sin\alpha  \theta^1 - a \cot\beta \theta^2 \\
\theta_4^5 & = & \phantom{-}  - a \cot\beta \csc\beta \theta^1 + \sin\alpha (\cot^2\beta -1)\theta^2,\nonumber
\end{eqnarray}
while the Gauss equation is equivalent to the equation:
\begin{equation}\label{eq: Gauss}
\begin{array}{lcl}
d\theta_2^1 + \theta_k^1 \wedge \theta_2^k & = & \theta^1 \wedge \theta^2.
\end{array}
\end{equation}
Therefore, using  equations  (\ref{eq: segunda}) and (\ref{eq:
Gauss}), we have
\begin{eqnarray}\label{eq:curvaturagauss}
K  & = &  1 - |\nabla \beta|^2 - 2 \cos \alpha \beta_1 - \cos^2 \alpha
-(1+ \csc^2 \beta )a^2  - \sin^2\alpha \cot^2\beta \\
&    =  &  1 - (1+ csc^2 \beta) a^2 \nonumber - |\nabla \beta + \cos \alpha e_1|^2 - \sin^2\alpha \cot\beta^2
\end{eqnarray}
Using (\ref{eq:intrin}) and  the complex structure of $S$, we get
\begin{equation} \label{eq:conex}
\begin{array}{lcl}
\theta_2^1  &  =  &  \tan\beta(d\beta\circ J-2\cos\alpha\theta^2)
\end{array}
\end{equation}
Differentiating (\ref{eq:conex}), we conclude that
\begin{eqnarray}
d\theta_2^1 & = & ( -(1 + \tan^2
\beta)|\nabla\beta|^2-\tan\beta\Delta\beta
-2\cos\alpha(1+2\tan^{2}\beta)\beta_1\nonumber\\
    &&  -4\tan^{2}\beta\cos^{2}\alpha) \theta^1
    \wedge \theta^2 \nonumber
\end{eqnarray}
where $\Delta = tr \nabla^2 $ is the Laplacian of $S$. The Gaussian curvature is therefore
given by:
\begin{eqnarray}\label{eq:curvatura2}
K & = &  -(1 + \tan^2 \beta)|\nabla\beta|^2-\tan\beta\Delta\beta
-2\cos\alpha(1+2\tan^{2}\beta)\beta_1 \nonumber\\
    &&
    -4\tan^{2}\beta\cos^{2}\alpha.
\end{eqnarray}
From  (\ref{eq:curvaturagauss}) and (\ref{eq:curvatura2}), we obtain
the following formula for the Laplacian of $S$:
\begin{eqnarray}\label{eq:lapla}
\tan\beta\Delta\beta &  = & (1+ \csc^2 \beta)a^2 \nonumber \\
&& - \tan^2 \beta( |\nabla \beta + 2 \cos \alpha e_1|^2
- |\sin \alpha ( 1- \cot^2 \beta)|^2 )\nonumber \\
 && + \sin^2 \alpha ( 1 -\tan^2 \beta)
\end{eqnarray}
\section{Proof of the Results}\label{eq:total}
In this section, in order to compute the Gauss-Codazzi-Ricci equations, we consider the holomorphic angle ($\alpha$) constant, and the principal curvature in the direction of $e_3$ constant, that is, $a$ is constant.
\subsection{Proof of Theorem\ref{teorema1}}
Now Codazzi-Ricci equations:
\begin{eqnarray}
d\theta_1^3 + \theta_2^3 \wedge \theta_1^2 + \theta_4^3 \wedge \theta_1^4 + \theta_5^3 \wedge \theta_1^5  &  =  & 0 \nonumber\\
d\theta_2^4 + \theta_1^4 \wedge \theta_2^1 + \theta_3^4 \wedge \theta_2^3 + \theta_5^4 \wedge \theta_2^5  &  =  & 0 \nonumber\\
d\theta_4^5 + \theta_1^5 \wedge \theta_4^1 + \theta_2^5 \wedge \theta_4^2 + \theta_3^5 \wedge \theta_4^3  &  =  & 0 \nonumber
\end{eqnarray}
simplify to:
\begin{eqnarray}\label{eq:projecao}
\beta_2 = \frac{(3-\cos^2\beta)a}{\sin\beta\cos\beta}(-2\sin\alpha\csc\beta\beta_1-\sin\alpha\cos\alpha\csc\beta(3-\cot^2\beta)+a^2\cot\alpha\csc\beta\cot^2\beta)
\end{eqnarray}
The following Codazzi-Ricci equations:
\begin{eqnarray}
d\theta_2^3 + \theta_1^3 \wedge \theta_2^1 + \theta_4^3 \wedge \theta_2^4 + \theta_5^3 \wedge \theta_2^5  &  =  & 0 \nonumber \\
d\theta_1^4 + \theta_2^4 \wedge \theta_1^2 + \theta_3^4 \wedge \theta_1^3 + \theta_5^4 \wedge \theta_1^5  &  =  & 0 \nonumber \\
d\theta_3^5 + \theta_1^5 \wedge \theta_3^1 + \theta_2^5 \wedge \theta_3^2 + \theta_4^5 \wedge \theta_3^4  &  =  & 0 \nonumber\\
d\theta_2^5 + \theta_1^5 \wedge \theta_2^1 + \theta_3^5 \wedge \theta_2^3 + \theta_4^5 \wedge \theta_2^4  &  =  & 0 \nonumber
\end{eqnarray}
reduces to:
\begin{eqnarray}\label{eq:projecao1}
\beta_1 = -  2\cos\alpha
\end{eqnarray}
Gauss-Codazzi-Ricci equations:
\begin{eqnarray}
d\theta_1^2 + \theta_3^2 \wedge \theta_1^3 + \theta_4^2 \wedge \theta_1^4 + \theta_5^2 \wedge \theta_1^5  &  =  & \theta^2 \wedge \theta^1 \nonumber \\
d\theta_1^5 + \theta_2^5 \wedge \theta_1^2 + \theta_3^5 \wedge \theta_1^3 + \theta_4^5 \wedge \theta_1^4  &  =  & 0 \nonumber\\
d\theta_3^4 + \theta_1^4 \wedge \theta_3^1 + \theta_2^4 \wedge \theta_3^2 + \theta_5^4 \wedge \theta_3^5 &  =  & 0 \nonumber
\end{eqnarray}
give the following Laplacian equation of the Contact angle ($\beta$):
\begin{eqnarray}\label{eq:lapangle}
\tan\beta\Delta\beta &  = & (1+ \csc^2 \beta)a^2 \nonumber \\
&& - \tan^2 \beta( |\nabla \beta + 2 \cos \alpha|^2
- |\sin \alpha ( 1- \cot^2 \beta)|^2 )\nonumber \\
 && + \sin^2 \alpha ( 1 -\tan^2 \beta)
\end{eqnarray}
Using (\ref{eq:projecao}), (\ref{eq:projecao1}) and (\ref{eq:lapangle}), we obtain the Theorem \ref{teorema1}. {$\square$} \\
\subsection{Proof of Corollary \ref{corolario1}}
Using (\ref{eq:projecao1}) in equation (\ref{eq:curvatura2}), we have:
\begin{eqnarray}
K = -(1+\tan^2\beta)\beta^2_2 - \tan\beta \Delta\beta
\end{eqnarray}
Therefore:
\begin{eqnarray}
\tan\beta \Delta\beta =  - K - (1+\tan^2\beta) \beta^2_2
\end{eqnarray}
Now using the condition that $K \geq 0$ and the Hopf's Lemma ( for $0 < \beta < \pi/2$), we get that the contact angle $\beta$ is constant, which prove the Corollary \ref{corolario1}. $\square$.
\begin{Remark}
The above result extends  the results in \cite{RMV}, in the sense that we don't need suppose that both angles are constant, and also the condition that the Gaussian Curvature is not constant, just non-negative.
\end{Remark}
{\bf ACKNOWLEDGEMENTS.} I would like to thank Prof. Gary R. Jensen for his valuable encouragement and suggestions during this work. I would also like to thank Department of Mathematics at Washington University in Saint Louis (WUSL) for their hospitality, and the Brazilian National Research Council (CNPq) for his support.


\end{document}